\documentclass[final,5p,times,twocolumn]{elsarticle}

\usepackage{amssymb}
\usepackage{amsmath}
\usepackage{graphics} % for pdf, bitmapped graphics files
\usepackage{epsfig} % for postscript graphics files
\usepackage{stfloats}
\usepackage{mathtools}
\usepackage{color}
\usepackage{multirow}
\usepackage{hhline}

\usepackage{lineno}

\journal{ArXiv}

\begin{document}

\begin{frontmatter}

%% Title, authors and addresses

%% use the tnoteref command within \title for footnotes;
%% use the tnotetext command for theassociated footnote;
%% use the fnref command within \author or \address for footnotes;
%% use the fntext command for theassociated footnote;
%% use the corref command within \author for corresponding author footnotes;
%% use the cortext command for theassociated footnote;
%% use the ead command for the email address,
%% and the form \ead[url] for the home page:
%% \title{Title\tnoteref{label1}}
%% \tnotetext[label1]{}
%% \author{Name\corref{cor1}\fnref{label2}}
%% \ead{email address}
%% \ead[url]{home page}
%% \fntext[label2]{}
%% \cortext[cor1]{}
%% \affiliation{organization={},
%%             addressline={},
%%             city={},
%%             postcode={},
%%             state={},
%%             country={}}
%% \fntext[label3]{}

\title{Simulation of Transients in Natural Gas Networks via A Semi-analytical Solution Approach}

%% use optional labels to link authors explicitly to addresses:
%% \author[label1,label2]{}
%% \affiliation[label1]{organization={},
%%             addressline={},
%%             city={},
%%             postcode={},
%%             state={},
%%             country={}}
%%
%% \affiliation[label2]{organization={},
%%             addressline={},
%%             city={},
%%             postcode={},
%%             state={},
%%             country={}}

\author[inst2,inst3]{Xin Xu\fnref{fn1}}
\author[inst2]{Rui Yao\corref{cor1}}
\author[inst3]{Kai Sun}
\author[inst2]{Feng Qiu}

\cortext[cor1]{Corresponding author}
\fntext[fn1]{This author is currently with Dominion Energy, Richmond, 23220, VA, USA.}

\affiliation[inst2]{organization={Energy System Division, Argonne National Laboratory},%Department and Organization
            % addressline={}, 
            city={Lemont},
            postcode={60439}, 
            state={IL},
            country={USA}}
            
\affiliation[inst3]{organization={Dept. of Electrical Engineering \& Computer Science, the University of Tennessee},%Department and Organization
            % addressline={}, 
            city={Knoxville},
            postcode={37996}, 
            state={TN},
            country={USA}}

\begin{abstract}
%% Text of abstract
Simulation and control of the transient flow in natural gas networks involve solving partial differential equations (PDEs). This paper proposes a semi-analytical solutions (SAS) approach for fast and accurate simulation of the natural gas transients. The region of interest is divided into a grid, and an SAS is derived for each grid cell in the form of the multivariate polynomials, of which the coefficients are identified according to the initial value and boundary value conditions. The solutions are solved in a ``time-stepping'' manner; that is, within one time step, the coefficients of the SAS are identified and the initial value of the next time step is evaluated. This approach achieves a much larger grid cell than the widely used finite difference method, and thus enhances the computational efficiency significantly. To further reduce the computation burden, the nonlinear terms in the model are simplified, which induces another SAS scheme that can greatly reduce the time consumption and have minor impact on accuracy. The simulation results on a single pipeline case and a 6-node network case validate the advantages of the proposed SAS approach in accuracy and computational efficiency.
\end{abstract}

% %%Graphical abstract
% \begin{graphicalabstract}
% %\includegraphics{grabs}
% \end{graphicalabstract}

%%Research highlights
\begin{highlights}
\item Research highlight 1: We present a semi-analytical solution approach for the simulation of transients in natural gas networks. 

\item Research highlight 2: To further reduce the computation burden, the nonlinear terms in the model are simplified which induces another SAS scheme that can greatly reduce the time consumption and have minor impact on accuracy.

\item Research highlight 3: We analyze the proposed approach in terms of efficiency and accuracy, and compare it with with the finite difference method.

\end{highlights}

\begin{keyword}
%% keywords here, in the form: keyword \sep keyword
Natural gas network \sep Semi-analytic approach \sep Simulation \sep Transient flow
%% PACS codes here, in the form: \PACS code \sep code
%% \PACS 0000 \sep 1111
%% MSC codes here, in the form: \MSC code \sep code
%% or \MSC[2008] code \sep code (2000 is the default)
%% \MSC 0000 \sep 1111
\end{keyword}

\end{frontmatter}

%% \linenumbers

%% main text
\section{Introduction}
%\label{sec:sample1}
Natural gas is an important energy source that can be used for electricity production and heating supply \cite{ma2019modeling, zlotnik2015optimal, chiang2016large}. As one of the crucial infrastructures for gas delivery, the gas pipeline networks have been investigated over decades on economic operation, control and security assessment. Disturbances in operating natural gas network can lead to transient processes propagating across the entire network, and potentially could cause damage to the equipment and even more severe security threats. Therefore, it is necessary to model and simulate the transient processes of natural gas pipeline network under potential disturbances. Simulations of the natural gas transient can directly reflect the variation of the states like pressure within the gas pipelines and provide important insight in the dynamics analysis and the system control design, which motivates numerous relevant studies \cite{stoner1969analysis, thorley1987unsteady, helgaker2014transient}. 

Natural gas pipeline transients can generally be formulated as a boundary value problem (BVP) of partial differential equations (PDEs) which governs the dynamics within the pipelines as well as network and terminal constraints. Many simulation techniques have been developed and investigated to solve the BVP. Two mainstream representatives are the method of characteristics \cite{stoner1969analysis} and the finite difference method (FDM) \cite{ thorley1987unsteady, helgaker2014transient, streeter1970natural, osiadacz1984simulation}. The method of characteristics is one of the early efforts to solve the BVP. Although an accurate result can be reached, it can be computationally intensive since tiny time steps are necessary to keep the numerical stability \cite{streeter1970natural}. On the other hand, the FDM gains more attention in recent years since it is more efficient and easier to implement \cite{thorley1987unsteady}. Other techniques are also investigated like finite volume method (FVM) \cite{qiu2018efficient} and state-space method \cite{alamian2012state},
which are still under development and need more comprehensive investigation to justify their performance in terms of accuracy and efficiency %\rcolor{The basic idea of FVM is to first convert the PDEs to ordinary differential equations (ODEs) via integral operations, and then, the result is obtained by solving the ODEs. Although \cite{qiu2018efficient} claims that FVM is more efficient than FDM, the benchmark FDM formulation in \cite{qiu2018efficient} is not presented which makes the comparison result in doubt. The state-space method is to convert the PDEs to transfer functions via linearization technique and Laplace transformation. However, this method does not show }

Amongst all those methods, the FDM is one of the most commonly used in both academia and industry. FDM is essentially a numerical method that depends on the discretization techniques, i.e. a grid should be selected to discretize the spatio-temporal space. The grid is commonly a regular one with square and rectangle cells. The basic idea of FDM is to approximate the partial derivatives in the PDEs by the values at the grid nodes, and then the BVP is converted to a set of difference equations, whose solution is the approximate solution of the BVP. The formulation of FDM can be either explicit or implicit in time. The explicit formulation can be easier to implement, but is usually limited to small time steps due to the stability confinements. The implicit formulation is more complex to implement but it has better numerical stability and thus enables larger time and spatial steps, which makes it the choice of many commercial software tools \cite{helgaker2013modeling}. Various studies have been conducted based on the FDM in terms of modeling\cite{herran2009modeling, helgaker2013modeling}, formulations \cite{thorley1987unsteady,helgaker2014transient}, and accuracy and numerical stability \cite{kiuchi1994implicit, greyvenstein2002implicit}. Note that one bottleneck for improving the efficiency of FDM lies in the grid cell size. Since all the constraints in the BVP are only enforced on the grid nodes, the grid cell size has to be limited in order to ensure the accuracy, which in turn  increases computation burden and limits the improvement of efficiency. 

This paper proposes a semi-analytical solution (SAS) approach for the simulation of transients in natural gas pipeline networks. The name ``semi-analytical'' comes from the algorithm that first sets a solution in a piece-wise analytical form (multivariate polynomials) and then determines the coefficients in the polynomials, which approximates the solution of the studied PDE. As will be shown later, SAS is more accurate and more efficient than FDM, because the constraints are enforced on the whole grid instead of just the gird nodes and also a larger grid cell can be achieved. In our previous works, SAS approach has been verified to be more efficient than many conventional numerical methods in solving ordinary differential equations (ODEs) \cite{8274026}, differential-algebraic equations (DAEs) \cite{wang2018time, yao2019efficient}, and partial differential equations \cite{xu2021semi}. The idea of SAS is also applied for solving algebraic equations in \cite{trias2012holomorphic, liu2017multi}, and \cite{xu2018holomorphic}. 

The rest of the paper is organized as follows. The BVP regarding the simulation of transient flow is formulated in section \ref{sec:2}. The algorithm of the SAS approach is introduced in section \ref{sec:3}. The overall simulation procedure is organized in section \ref{sec:4}. The simulation results on a single pipeline case and a 6-node network case are presented in section \ref{sec:5} and compared with those of FDM. The conclusions and future works are in section \ref{sec:6}.

\section{Boundary value problem of natural gas network transients}\label{sec:2}

\subsection{Natural gas network model}

The natural gas network can be viewed as a graph with links and nodes. The set of links $\mathcal{E}$ includes all the pipelines and the set of nodes $\mathcal{V}$ is the union of the following sets: (\textit{i}) supply nodes $\mathcal{V_P}$ where the gas is supplied into the network and the pressure $p$ is usually controlled, (\textit{ii}) demand nodes $\mathcal{V_Q}$ where the gas is extracted out of the network and the mass flow $q$ is determined, and (\textit{iii}) $\mathcal{V_J}$ junction nodes  that are not in $\mathcal{V_P}$ or $\mathcal{V_Q}$ and where pipelines are connected. 

For the sake of convenience, when introducing the network model, we use the superscripts $(e)$ and $(\nu)$ to denote the associated quantity/function of the pipeline $e \in \mathcal{E}$ and node $\nu \in \mathcal{V}$, respectively. 

The simplified hydraulic model of the pipeline is considered which governs the transient flow within the pipeline \cite{ helgaker2014transient, helgaker2013modeling}. Assume that pipelines are horizontal. For a segment of pipeline $e \in \mathcal{E}$ of the length $L^{(e)}$, the transients can be modeled by partial differential equations (PDEs) \eqref{eq_PDE}. All the physical quantities and the units are defined in Table \ref{table_Gas_quan}. Without loss of generality, we choose the convention that the gas flows from $x=0$ to $x=L^{(e)}$. 

\begin{equation} \label{eq_PDE}
    \begin{cases}
        \partial_t{p^{(e)}} + \frac{v^2}{S^{(e)}}\partial_x{q^{(e)}} = 0\\
        \partial_t{q^{(e)}} + S^{(e)} \partial_x{p^{(e)}} + \frac{\lambda^{(e)} v^2 q^{(e)} |q^{(e)}|}{2 d^{(e)} S^{(e)} p^{(e)}} = 0
    \end{cases}
\end{equation}

\begin{table}
\caption{Physical Quantities of Gas Pipeline}
\label{table_Gas_quan}
\begin{center}
\begin{tabular}{|c||c|}
\hline
Physical quantity & Units\\
\hline
Pressure, $p$ & $\mathrm{Pa}$\\
\hline
Mass flow, $q$ &  $\mathrm{kg/s}$\\
\hline
Sound speed, $\nu$ & $\mathrm{m/s}$ \\
\hline
Cross-section area of pipeline, $S$ &  $\mathrm{m}^2$\\
\hline
Diameter of pipeline, $d$ &  $\mathrm{m}$\\
\hline
Friction factor, $\lambda$ & - \\
\hline
Constant temperature, $T_0$ &  $\mathrm{K}$\\
\hline
Specific gas constant, $R$ &  $\mathrm{J/(kg \cdot K)}$\\
\hline
\end{tabular}
\end{center}
\end{table}

The gas transients are also determined by the given initial value and the constraints imposed at the nodes. The initial value of the pipeline $e \in \mathcal{E}$ is given at $t=0$ for the any location $x \in [0, L^{(e)}]$, as shown in \eqref{eq_ini} where $P^{(e)}_{ini}$ and $Q^{(e)}_{ini}$ defines the initial value.

\begin{equation} \label{eq_ini}
    \begin{cases}
        p^{(e)}(x,0) = P^{(e)}_{ini}(x)\\
        q^{(e)}(x,0) = Q^{(e)}_{ini}(x)
    \end{cases}
\end{equation}

The constraints imposed at a node $\nu \in \mathcal{V}$ should consider both the specific node and also the associated pipelines. For a pipeline $e \in \mathcal{E}$ whose inlet is connected to a supplying node $\nu \in \mathcal{V_P}$, the controlled $p$ follows a Dirichlet boundary condition as in \eqref{eq_bd_P}, where $P^{(\nu)}_B$ defines the boundary condition. 

\begin{equation} \label{eq_bd_P}
        p^{(e)}(0,t) = P^{(\nu)}_B(t) 
\end{equation}

If the outlet of the pipeline $e \in \mathcal{E}$ is connected to a demanding node $\nu \in \mathcal{V_Q}$, the controlled $q$ follows a Dirichlet boundary condition as in \eqref{eq_bd_Q}, where $Q^{(\nu)}_B$ defines the boundary condition.

\begin{equation} \label{eq_bd_Q}
    q^{(e)}(L^{(e)},t) = Q^{(\nu)}_B(t)
\end{equation}

If the pipelines $e_{in,1}, e_{in,2}, ...$ have their inlet connected to a junction node $v \in \mathcal{V_J}$ and the pipeline $e_{out,1}, e_{out,2}, ...$ have their outlet connected to the same junction node, \eqref{eq_bd_J} must be satisfied, i.e. at the junction node, the pressure of all the pipelines should be the same, and the mass flow should be balanced. $Q^{(\nu)}_J$ determines the amount of gas extracted out of the network at node $v$.

\begin{equation} \label{eq_bd_J}
    \begin{cases}
        p^{(e_{in,1})}(0,t) = p^{(e_{in,i})}(0,t),\;\;\;\;\;\;\;\;\;\;\;\;\;\;\;\;\;\;\;\;\;i \neq 1 \\
        p^{(e_{in,1})}(0,t) = p^{(e_{out,j})}(L^{(e_{out,j})},t)\\
        \sum_{j} q^{ (e_{out,j}) }(L^{ (e_{out,j}) },t) - \sum_{i} q^{(e_{in,i})}(0,t) = Q^{(\nu)}_J(t)
    \end{cases}
\end{equation}

\subsection{Grid selection and normalized BVP}

The boundary value problem of natural gas network is to identify the solution of pressure $p^{(e)}(x,t)$ and mass flow $q^{(e)}(x,t)$ regarding the temporal variable $t$ and spatial variable $x$, by considering the constraints below: 

\begin{itemize} 
    \item Transient flow in the pipeline governed by the PDEs, \eqref{eq_PDE}.
    \item Given initial values, \eqref{eq_ini}.
    \item Controlled $p$ at the supplying nodes $\mathcal{V_P}$, \eqref{eq_bd_P}.
    \item Controlled $q$ at the demanding nodes $\mathcal{V_Q}$, \eqref{eq_bd_Q}.
    \item Constraints at the junction nodes $\mathcal{V_J}$, \eqref{eq_bd_J}.
\end{itemize}

\begin{figure}
	\centering
		\includegraphics[scale=.5]{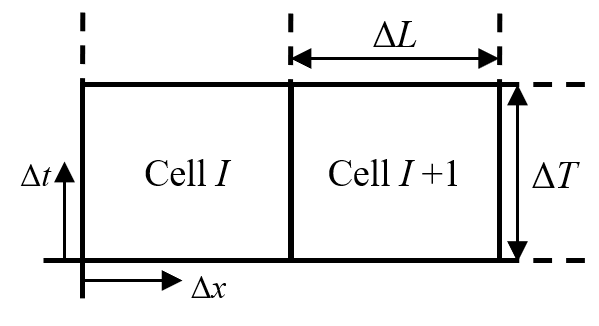}
	\caption{Illustration of grid cells in $x$-$t$ space.}
	\label{Fig1_grid}
\end{figure}

For each pipeline, the region  of interest is divided into a grid in the $x$-$t$ space like Fig. \ref{Fig1_grid}, which is a common practice for many numerical methods. Moreover, a uniform grid is selected, i.e. each cell has the same size $\Delta L \times \Delta T$. Assume that within one time step, there are $N$ uniform cells of size $\Delta L \times \Delta T$. Note that for different pipeline segments , say $(e_i)$ and $(e_j)$, it is required that $\Delta T^{(e_i)} = \Delta T^{(e_j)}$, i.e. the time steps are synchronized, but $\Delta L^{(e_i)} = \Delta L^{(e_j)}$ and $ N^{(e_i)} = N^{(e_j)}$ are not required. Henceforth, we only use $\Delta T$ to denote the temporal length without specifying the pipeline.

After defining the grid, the solutions $p^{(e,I)}(x,t)$ and $q^{(e,I)}(x,t)$ of each cell will constitute the solution of the entire pipeline, where $(e,I)$ denotes the Cell $I$ of the pipeline $e$. And thus for any two adjacent cells of each pipeline, their solutions should be continuous at the seam, which results in an additional "seamless" constraint.

\begin{equation} \label{eq_seam}
    \begin{cases}
        p^{(e,I)}(\Delta L^{(e)},\Delta t) = p^{\left(e,(I+1)\right)}(0,\Delta t)\\
        q^{(e,I)}(\Delta L^{(e)}, \Delta t) = q^{\left(e,(I+1)\right)}(0,\Delta t)
    \end{cases}
\end{equation}

Normalization technique is used to simplify the formulation and the follow-up analysis. Two base values $p_b$ and $q_b$ are selected to normalize the $p^{(e)}$ and $q^{(e)}$ of each pipeline $e$, respectively. Then, for the cells of each pipeline, $\Delta x$ and $\Delta t$ are normalized by $\Delta L^{(e)}$ and $\Delta T$, respectively, which transforms a grid cell to $[0, 1] \times [0, 1]$. The superscript ``*'' denotes the nominal values.

After normalization, \eqref{eq_PDE}, \eqref{eq_ini}, \eqref{eq_bd_P}, \eqref{eq_bd_Q}, \eqref{eq_bd_J} and \eqref{eq_seam} are normalized to \eqref{eq_nor_PDE}, \eqref{eq_nor_ini}, \eqref{eq_nor_bd_P}, \eqref{eq_nor_bd_Q}, \eqref{eq_nor_bd_J} and \eqref{eq_nor_seam}, respectively. Note that the superscript ``*'' is omitted to simplify the expression. $P^{(e,I)}_{ini,t}$ is derived from $P^{(e)}_{ini}$ when the starting time of the Cell $I$ is $t$. The rest of the functions with the subscript $t$ are also similarly derived. 

\begin{equation} \label{eq_nor_PDE}
    \begin{cases}
        \partial_{\Delta t} p^{(e,I)} + C^{(e)}_1 \partial_{\Delta x} q^{(e,I)} = 0 \\
        \partial_{\Delta t} q^{(e,I)} + C^{(e)}_2 \partial_{\Delta x} p^{(e,I)} + C^{(e)}_3 \frac{ q^{(e,I)} |q^{(e,I)}| }{p^{(e,I)}} = 0
    \end{cases}
\end{equation}

\begin{equation} \label{eq_nor_C}
\begin{cases}
    C^{(e)}_1 = \frac{v^2 q_b \Delta T}{S^{(e)} p_b \Delta L^{(e)}}\\ C^{(e)}_2 = \frac{S^{(e)} p_b \Delta T}{q_b \Delta L^{(e)}} \\
    C^{(e)}_3 = \frac{\lambda^{(e)} v^2 q_b \Delta T}{2 d^{(e)} S^{(e)} p_b }
\end{cases}
\end{equation}

\begin{equation} \label{eq_nor_ini}
    \begin{cases}
        p^{(e,I)}(\Delta x,0) = P^{(e,I)}_{ini,t}(\Delta x)\\
        q^{(e,I)}(\Delta x,0) = Q^{(e,I)}_{ini,t}(\Delta x)
    \end{cases}
\end{equation}

\begin{equation} \label{eq_nor_bd_P}
    p^{(e,1)}(0,\Delta t) = P^{(\nu)}_{B,t}(\Delta t)
\end{equation}

\begin{equation} \label{eq_nor_bd_Q}
    q^{\left(e,N^{(e)}\right)}(1, \Delta t) = Q^{(\nu)}_{B,t}(\Delta t)
\end{equation}

\begin{equation} \label{eq_nor_bd_J}
    \begin{cases}
        p^{(e_{in,1},1)}(0,\Delta t) = p^{(e_{in,i},1)}(0,\Delta t),\;\;\;\;\;\;\;\;\;\;\;\;\;\;i \neq 1 \\
        p^{(e_{in,1},1)}(0,\Delta t) = p^{\left(e_{out,j}, N^{(e_{out,j})}\right)}(1,\Delta t)\\
        \sum_{e_{out,j}} q^{ \left(e_{out,j},N^{(e_{out,j})}\right) }( 1, \Delta t) - \\
        \;\;\;\;\;\;\;\;\;\;\;\;\;\;\;\;\;\;\;\;\;\;\;\sum_{e_{in,i}} q^{(e_{in,i},1)}(0,\Delta t) = Q^{(\nu)}_{J,t}(\Delta t) 
    \end{cases}
\end{equation}

\begin{equation} \label{eq_nor_seam}
    \begin{cases}
        p^{(e,I)}(1, \Delta t) = p^{\left(e,(I+1)\right)}(0,\Delta t) \\
        q^{(e,I)}(1, \Delta t) = q^{\left(e,(I+1)\right)}(0,\Delta t)
    \end{cases}
\end{equation}

Henceforth, the original BVP is converted to finding the solution of all the cells under the constraints of \eqref{eq_nor_PDE}, \eqref{eq_nor_ini}, \eqref{eq_nor_bd_P}, \eqref{eq_nor_bd_Q}, \eqref{eq_nor_bd_J} and \eqref{eq_nor_seam}.

\section{Semi-analytical solution approach}\label{sec:3}

\subsection{Embedding technique and semi-analytical solution}

Without loss of generality, only consider the case $q^{(e,I)} \geq 0$. Hence, the absolute operator $|\cdot|$ in \eqref{eq_nor_PDE} can be ignored. Since \eqref{eq_nor_PDE} is nonlinear, the embedding technique is used to reformulate it such that the SAS approach can be easily applied, as shown in \eqref{eq_ebd_PDE}.

\begin{equation} \label{eq_ebd_PDE} \small
    \begin{cases}
        \partial_{\Delta t} p^{(e,I)} + C^{(e)}_1 \partial_{\Delta x} q^{(e,I)} = 0 \\
        p^{(e,I)} \partial_{\Delta t} q^{(e,I)} + C^{(e)}_2 p^{(e,I)} \partial_{\Delta x} p^{(e,I)} + s C^{(e)}_3 \left(q^{(e,I)}\right)^2 = 0
    \end{cases}
\end{equation}

\noindent where an embedding variable $s$ is inserted at the friction term i.e. the term regarding $C^{(e)}_3$. $s$ can be viewed as a scaling factor of the friction factor $\lambda^{(e)}$. When $s=0$, the friction vanishes and \eqref{eq_ebd_PDE} can actually be reduced to constant coefficient first-order PDEs \eqref{eq_ebd_PDE_s0}. On the other hand, when $s = 1$, \eqref{eq_ebd_PDE} is equivalent to the original PDEs. 

\begin{equation} \label{eq_ebd_PDE_s0}
    \begin{cases}
        \partial_{\Delta t} p^{(e,I)} + C^{(e)}_1 \partial_{\Delta x} q^{(e,I)} = 0 \\
        \partial_{\Delta t} q^{(e,I)} + C^{(e)}_2 \partial_{\Delta x} p^{(e,I)} = 0
    \end{cases}
\end{equation}

The embedding variable $s$ enables more flexibility in dealing with the nonlinearity in the original PDEs; that is, instead of directly solving the original problem, we can first ignore the nonlinear friction term by setting $s=0$ and solve the simplified problem, and then, based the solution at $s=0$, we can reconsider the friction factor by considering $s>0$ and identify the solution of the original PDEs at $s=1$. 

\subsection{Methodology of SAS approach}

After embedding, the solution of each cell, $p^{(e,I)}$ and $q^{(e,I)}$, is approximated by SAS, i.e. a multivariate polynomial regrading $\Delta x$, $\Delta t$ and $s$, as defined in \eqref{eq_SAS_Sol}. Here, $R$ represents the maximum order of $s^r$ and, when $r$ is fixed, $M$ represents the maximum order of the polynomials regarding $\Delta x^n \Delta t^{m-n}$.

\begin{equation} \label{eq_SAS_Sol}
    \begin{aligned}
        p^{(e,I)}(\Delta x, \Delta t) = \sum^{R}_{r=0} \left[ \sum^{M}_{m=0}\sum^{m}_{n = 0} p^{(e,I)}_{n,m-n,r}\Delta x^n \Delta t^{m-n}  \right] s^r\\
        q^{(e,I)}(\Delta x, \Delta t) = \sum^{R}_{r=0} \left[ \sum^{M}_{m=0}\sum^{m}_{n = 0} q^{(e,I)}_{n,m-n,r}\Delta x^n \Delta t^{m-n}  \right] s^r
    \end{aligned}
\end{equation}

\noindent where $p^{(e,I)}_{i,j,r}$ and $q^{(e,I)}_{i,j,r}$ are scalar coefficients and also the unknowns to be identified. To reduce the computational complexity, some restrictions are applied on the terms $p^{(e,I)}_{0,0,r}$ and $q^{(e,I)}_{0,0,r}$ such that they can be treated as given values instead of unknowns. Specifically, set $p^{(e,I)}_{0,0,0} = p^{(e,I)}_{ini,t}(0)$ and $q^{(e,I)}_{0,0,0} = q^{(e,I)}_{ini,t}(0)$ which are immediately given by the initial value, and set $p^{(e,I)}_{0,0,r} = q^{(e,I)}_{0,0,r} = 0$ for $r>0$. By doing so, the solution of \eqref{eq_ebd_PDE} at the initial time is always equal to the given initial value regardless of $s= 0$ or 1. The rest of the coefficients $p^{(e,I)}_{i,j,r}$ and $q^{(e,I)}_{i,j,r}$ need to be calculated so that the solutions \eqref{eq_SAS_Sol} satisfy the constraints of BVP. For the sake of convenience, define that $p^{(e,I)}_{i,j,r}$ and $q^{(e,I)}_{i,j,r}$ are the unknowns regarding $s^r$. When $r$ is fixed, there are totally $M(M+3)$ unknowns regarding $s^r$ for each cell. The total number of the unknowns are $(R+1)M(M+3)$ for each cell. 

The SAS approach derives the solution in a ``time-stepping'' manner, i.e. along the temporal axis, it will iteratively solve the solution of each row of grid cells from the start time to the designated end time. The initial value of each row will be given by evaluating the semi-analytical solutions of the previous row at $\Delta t = 1$.

Within each time step, the unknowns will be identified in a recursive manner. The unknowns regarding $s^0$ are firstly identified, and then, those regarding $s^1$, $s^2$ and so on, i.e. identify the unknowns from $s^0 \rightarrow s^1 \rightarrow s^1 \dots \rightarrow s^R$.

The subsequent subsections will discuss the details of the recursive solution of the unknowns within each time step. As aforementioned, we will first set $s = 0$ to ignore the friction and solve the remaining problem, which is equivalent to solve all the unknowns regarding $s^0$. Then, the rest of unknowns can be recursively identified for each $s^r$. % At the end, how to reduce the error during the ``time-stepping'' simulation is discussed.

\subsection{Solution of unknowns regarding $s^r$, $r = 0$}

The simplified problem at $s = 0$ considers the constraints \eqref{eq_ebd_PDE_s0}, \eqref{eq_nor_ini}, \eqref{eq_nor_bd_P}, \eqref{eq_nor_bd_Q}, \eqref{eq_nor_bd_J} and \eqref{eq_nor_seam}. The solution \eqref{eq_SAS_Sol} is simplified to \eqref{eq_SAS_Sol_s0}, which shows that solving the simplified problem is equivalent to identifying the unknowns regarding $s^0$, $p^{(e,I)}_{n,m-n,0}$ and $q^{(e,I)}_{n,m-n,0}$, for the original problem. By substituting solution \eqref{eq_SAS_Sol_s0} into the constraints, a set of linear equations regarding the unknowns will be obtained and the unknowns can be solved. The treatments of the constraint for each cell is discussed as follows.

\begin{equation} \label{eq_SAS_Sol_s0}
    \begin{aligned}
        p^{(e,I)}(\Delta x, \Delta t) = \sum^{M}_{m=0}\sum^{m}_{n = 0} p^{(e,I)}_{n,m-n,0}\Delta x^n \Delta t^{m-n}  \\
        q^{(e,I)}(\Delta x, \Delta t) = \sum^{M}_{m=0}\sum^{m}_{n = 0} q^{(e,I)}_{n,m-n,0}\Delta x^n \Delta t^{m-n}
    \end{aligned}
\end{equation}

\subsubsection{PDE constraints}

In each cell, by substituting \eqref{eq_SAS_Sol_s0} into \eqref{eq_ebd_PDE_s0} and equating the coefficients of the resulting PDEs we can get:

\begin{equation} \label{eq_SAS_PDE_der1}
\begin{cases}
    \sum^{M}_{m=0}\sum^{m-1}_{n = 0} (m-n) p^{(e,I)}_{n,m-n,0}\Delta x^n \Delta t^{m-n-1}  + \\
    \;\;\;\;\;\;\;\;\;\;\;\;\;\;\;\; C^{(e)}_1 \sum^{M}_{m=0}\sum^{m}_{n = 1} n q^{(e,I)}_{n,m-n,0} \Delta x^{n-1} \Delta t^{m-n}  = 0 \\
    \sum^{M}_{m=0}\sum^{m-1}_{n = 0} (m-n) q^{(e,I)}_{n,m-n,0}\Delta x^n \Delta t^{m-n-1}  + \\
    \;\;\;\;\;\;\;\;\;\;\;\;\;\;\;\; C^{(e)}_2 \sum^{M}_{m=0}\sum^{m}_{n = 1} n p^{(e,I)}_{n,m-n,0} \Delta x^{n-1} \Delta t^{m-n}  = 0 \\
\end{cases}
\end{equation}
Then organize the polynomial terms with the same orders,

\begin{equation} \label{eq_SAS_PDE_der2}
\begin{cases}
    \sum^{M}_{m=0}\sum^{m-1}_{n = 0} \Big[ (m-n) p^{(e,I)}_{n,m-n,0}  + \\
    \;\;\;\;\;\;\;\;\;\;\;\;\;\;\;\; (n+1) C^{(e)}_1 q^{(e,I)}_{n+1,m-n-1,0} \Big]  \Delta x^n \Delta t^{m-n-1}  = 0 \\
    \sum^{M}_{m=0}\sum^{m-1}_{n = 0} \Big[ (m-n) q^{(e,I)}_{n,m-n,0}  + \\
    \;\;\;\;\;\;\;\;\;\;\;\;\;\;\;\; (n+1) C^{(e)}_2 p^{(e,I)}_{n+1,m-n-1,0} \Big] \Delta x^n \Delta t^{m-n-1}  = 0 \\
\end{cases}
\end{equation}
For the same polynomial term regarding $\Delta x$ and $\Delta t$, a set of linear equations regarding the unknowns $p^{(e,I)}_{n,m-n,0}$ and $q^{(e,I)}_{n,m-n,0}$ can be obtained:

\begin{equation} \label{eq_SAS_PDE_coef}
\begin{gathered}
    \begin{cases}
        (m-n) p^{(e,I)}_{n, m-n, 0} + (n+1) C^{(e)}_1 q^{(e,I)}_{n+1, m-n-1, 0} = 0 \\
        (m-n) q^{(e,I)}_{n, m-n, 0} + (n+1) C^{(e)}_2 p^{(e,I)}_{n+1, m-n-1, 0} = 0 
    \end{cases} \\
    0\leq n \leq m-1,\; for\;\;\; 1 \leq m \leq M
\end{gathered}
\end{equation}

Therefore, following the procedure of \eqref{eq_SAS_PDE_der1}, \eqref{eq_SAS_PDE_der2} and \eqref{eq_SAS_PDE_coef}, for each cell, the resulting equation regarding each polynomial term $\Delta x^n \Delta t^{m-n-1}$ is given in \eqref{eq_SAS_PDE_coef}, where $0 \leq n \leq m-1$ and $1 \leq m \leq M$. 

For a pipeline $e \in \mathcal{E}$, \eqref{eq_SAS_PDE_coef} corresponds to $N^{(e)}M(M+1)$ equations. Note that \eqref{eq_SAS_PDE_coef} is enforced on the whole spatio-temporal space instead of just the grid nodes.

\subsubsection{Initial value constraints}

\begin{figure}[thpb]
    \centering
    \includegraphics[scale=0.5]{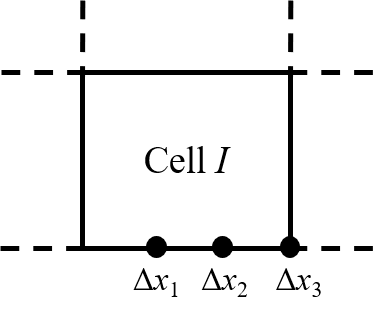}
    \caption{Selection of points for initial value constraints.}
    \label{fig:ini}
\end{figure}

Instead of requiring \eqref{eq_nor_ini} to be strictly met along the entire boundary, we only consider several points to meet \eqref{eq_nor_ini} as illustrated in Fig. \ref{fig:ini}, with each point corresponding to a distinct location $\Delta x_k \neq 0$ (because $\Delta x_k = 0$ already has $p^{(e,I)}_{0,0,0} = P^{(e,I)}_{ini,t}(0)$ and $q^{(e,I)}_{0,0,0} = Q^{(e,I)}_{ini,t}(0)$):

\begin{equation} \label{eq_SAS_ini_coef}
    \begin{cases}
        \sum^{M}_{m=0} \Delta x_k^m p^{(e,I)}_{m,0,0} = P^{(e,I)}_{ini,t} (\Delta x_k)\\
        \sum^{M}_{m=0} \Delta x_k^m q^{(e,I)}_{m,0,0} = Q^{(e,I)}_{ini,t} (\Delta x_k)  
    \end{cases} 
\end{equation}

For a pipeline $e \in \mathcal{E}$, if there are $M^{(e)}_x$ such points on each cell, \eqref{eq_SAS_ini_coef} corresponds to $N^{(e)}M^{(e)}_x$ linear equations in total.

\subsubsection{Boundary condition constraints regarding $\mathcal{V_P}$ and $\mathcal{V_Q}$}

\begin{figure}[thpb]
    \centering
    \includegraphics[scale=0.5]{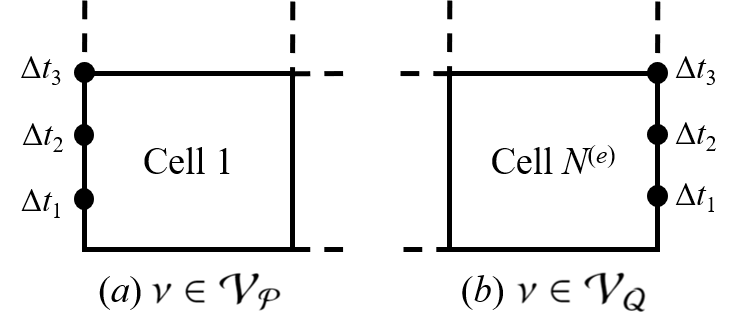}
    \caption{Selection of points for boundary condition constraints regarding $\mathcal{V_P}$ and $\mathcal{V_Q}$.}
    \label{fig:bon_P_Q}
\end{figure}

The boundary condition regarding a node $\nu \in \mathcal{V_P}$ is imposed on the Cell $1$ of the connected pipeline at $\Delta x = 0$,  and those regarding a node $\nu \in \mathcal{V_Q}$ is imposed on the Cell $N^{(e)}$ of the connected pipeline at $\Delta x = 1$. Similar to initial value constraint, only consider several points to meet \eqref{eq_nor_bd_P} and \eqref{eq_nor_bd_Q} as illustrated in Fig. \ref{fig:bon_P_Q}, with each point corresponding to a distinct time instant $\Delta t_k \neq 0$:

\begin{equation} \label{eq_SAS_bd_P_Q_coef}
    \begin{cases}
        \sum^{M}_{m=0} \Delta t_k^{m} p^{(e,1)}_{0,m,0}  = P^{(\nu)}_{B,t}(\Delta t_k), \;\;\;\;\;\;\;\;\;\;\;\;\;\;\;\;\; \text{if $\nu \in \mathcal{V_P}$} \\
        \sum^{M}_{m=0}\sum^{m}_{n = 0} \Delta t_k^{m-n} q^{(e,N^{(e)})}_{n,m-n,0}  = Q^{(\nu)}_{B,t}(\Delta t_k),\;\;\;\;\text{if $\nu \in \mathcal{V_Q}$}
    \end{cases} 
\end{equation}

For a node $\nu \in \mathcal{V_P} \cup \mathcal{V_Q}$, if there are $M^{(\nu)}_B$ such points on the boundary, \eqref{eq_SAS_bd_P_Q_coef}  corresponds to $M^{(\nu)}_B$ linear equations.

\subsubsection{Boundary condition constraints regarding $\mathcal{V_J}$}

The boundary condition regarding a node $\nu \in \mathcal{V_J}$ is imposed on the cells of multiple connected pipelines. For the pipelines $e_{in,1}, e_{in,2}, ...$ that have their inlet connected to $v$, the boundary condition is imposed to the Cell $I=1$ at $\Delta x = 0$. For the pipeline $e_{out,1}, e_{out,2}, ...$ have their outlet connected to $\nu$, the boundary condition is imposed on the Cell $I=N^{(e_{out,i})}$ at $\Delta x = 1$. Similar to the boundary condition constraints regarding $\mathcal{V_P}$ and $\mathcal{V_Q}$, only consider several points to meet \eqref{eq_nor_bd_J}, with each point corresponding to a distinct time instant $\Delta t_k \neq 0$:

\begin{equation} \label{eq_SAS_bd_J_coef}
\begin{gathered}
    \begin{cases}
        \sum^{M}_{m=0} \Delta t_k^m  p^{(e^{in,1},1)}_{0,m,0} = \sum^{M}_{m=0} \Delta t_k^m  p^{(e^{in,i},1)}_{0,m,0},\;\;\;\;\;\;\;\;\;\;\;i \neq 1 \\
        \sum^{M}_{m=0} \Delta t_k^m  p^{(e^{in,1},1)}_{0,m,0} = \sum^{M}_{m=0} \sum^{m}_{n=0} \Delta t_k^{m-n} p^{\left(e^{out,j},N^{(e^{out,j})}\right)}_{n,m-n,0}  \\
        \sum_{e^{out,j}} \left[   \sum^{M}_{m=0} \sum^{m}_{n=0} \Delta t_k^{m-n} q^{\left(e^{out,j},N^{(e^{out,j})}\right)}_{n,m-n,0}   \right] - \\
        \;\;\;\;\;\;\;\;\;\;\;\;\;\;\;\;\;\;\;\;\sum_{e^{in,i}} \left[  \sum^{M}_{m=0} \Delta t_k^m  q^{(e^{in,i},1)}_{0,m,0} \right] = Q^{(\nu)}_{J,t}(\Delta t_k) 
    \end{cases}
\end{gathered}
\end{equation}

For a node $\nu \in \mathcal{V_J}$, if there are $M^{(\nu)}_J$ such points on the boundary, \eqref{eq_SAS_bd_J_coef}  corresponds to $M^{(\nu)}_J (N^{(\nu)}_{in} + N^{(\nu)}_{out})$ linear equations, where $N^{(\nu)}_{in}$ and $N^{(\nu)}_{out}$ are the total number of $e_{in,i}$ and $e_{out,j}$, respectively.

\subsubsection{Seamless Condition Constraints}

\begin{figure}[thpb]
    \centering
    \includegraphics[scale=0.5]{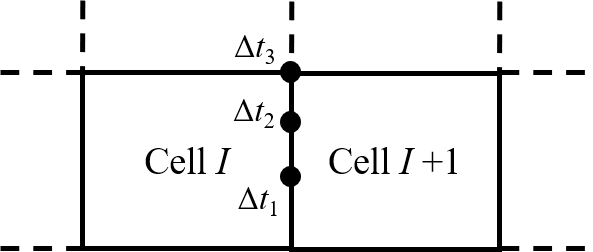}
    \caption{Selection of points for seamless condition constraints.}
    \label{fig:seam}
\end{figure}

Similar to the boundary condition constraints, only consider a limited number of points on the borders of the neighboring cells to meet \eqref{eq_nor_seam} as illustrated in Fig. \ref{fig:seam}, with each point corresponding to a distinct time instant $\Delta t_k \neq 0$:

\begin{equation} \label{eq_SAS_seam_coef}
\begin{gathered} 
    \begin{cases}
        \sum^{M}_{m=0}\sum^{m}_{n = 0} \Delta t_k^{m-n} p^{(e,I)}_{n,m-n,0} = \sum^{M}_{m=0} \Delta t_k^{m} p^{\left( e,(I+1) \right)}_{0,m,0}  \\
        \sum^{M}_{m=0}\sum^{m}_{n = 0} \Delta t_k^{m-n} q^{(e,I)}_{n,m-n,0} = \sum^{M}_{m=0} \Delta t_k^{m} q^{\left( e,(I+1) \right)}_{0,m,0}
    \end{cases}
\end{gathered}   
\end{equation}

For a pipeline $e \in \mathcal{E}$, if there are $M^{(e)}_s$ such points on the borders, \eqref{eq_SAS_seam_coef} corresponds to $(N^{(e)}-1) M^{(e)}_s$ linear equations in total.

\subsubsection{Linear Equations of SAS Coefficients}

All the linear equations from the constraints can be organized in the form of \eqref{eq_lin_sys}, where $\mathbf{c}$ is the vector of all the unknowns, $\mathbf{A}$ is the resulting matrix and $\mathbf{b}$ is a vector of constants.

\begin{equation} \label{eq_lin_sys}
    \mathbf{Ac} = \mathbf{b}
\end{equation}

The total number of equations and unknowns are computed by \eqref{eq_eqn_num} and \eqref{eq_unk_num}, respectively.

\begin{equation} \label{eq_eqn_num} \small
\begin{gathered} 
    N_{eqn} =  \sum_{e \in \mathcal{E}} \left( N^{(e)}M(M+1) + N^{(e)} M^{(e)}_x + (N^{(e)}-1)M^{(e)}_s \right) + \\
    \sum_{\nu \in \mathcal{V_P \cup V_Q}} M^{(\nu)}_B + \sum_{\nu \in \mathcal{V_J}} M^{(\nu)}_J \left(  N^{(\nu)}_{in} + N^{(\nu)}_{out}  \right) 
\end{gathered} 
\end{equation}

\begin{equation} \label{eq_unk_num}
    N_{unk} =  \sum_{e \in \mathcal{E}} \left( N^{(e)}M(M+3) \right) 
\end{equation}

\begin{figure*}
\begin{equation} \label{eq_SAS_PDE_coef_sr}
\begin{gathered}
    \begin{cases}
         (m-n) p^{(e,I)}_{n, m-n, r}  + (n+1) C^{(e)}_1 q^{(e,I)}_{n+1, m-n-1, r} = 0 \\
         \\
         \begin{aligned}
          (m-n) p^{(e,I)}_{0,0,0} q^{(e,I)}_{n, m-n, r} & + (n+1) C^{(e)}_2 p^{(e,I)}_{0,0,0} p^{(e,I)}_{n+1, m-n-1, r} + p^{(e,I)}_{n,m-n-1,0} q^{(e,I)}_{0,1,r} \left[  q^{(e,I)}_{0,1,0} + C^{(e)}_2  p^{(e,I)}_{1,0,0} \right]  p^{(e,I)}_{n,m-n-1,r} + C^{(e)}_2 p^{(e,I)}_{n,m-n-1,0}  p^{(e,I)}_{1,0,r} \\
          & +\sum^{m-2}_{\gamma=1} \sum^{\min(n,\gamma)}_{i=\max(0,n+\gamma-m+1)} \Bigg\{ \left[ (m-\gamma-n+i) q^{(e,I)}_{n-i,m-\gamma-n+i,0} + (n-i+1) C^{(e)}_2 p^{(e,I)}_{n-i+1, m-\gamma -n+i-1,0}  \right]  p^{(e,I)}_{i, \gamma-i,r} \\
         & + (m-\gamma-n+i) p^{(e,I)}_{i,\gamma-i,0} q^{(e,I)}_{n-i, m-\gamma-n+i,r}  + (n-i+1) C^{(e)}_2 p^{(e-I)}_{i,\gamma-i,0} p^{(e,I)}_{n-i+1, m-\gamma-n+i-1,r}  \Bigg\} \\
        = & -\sum^{r-1}_{j=1} p^{(e,I)}_{n,m-n-1,j}  q^{(e,I)}_{0,1,r-j} - C^{(e)}_2 \sum^{r-1}_{j=1} p^{(e,I)}_{n,m-n-1,j} p^{(e,I)}_{1,0,r-j} - 2 C^{(e)}_3 q^{(e,I)}_{0,0,0} q^{(e,I)}_{n,m-m-1,r-1}\\
        & -  \sum^{m-2}_{\gamma=1} \sum^{\min(n,\gamma)}_{i=\max(0,n+\gamma-m+1)} \Bigg\{ (m-\gamma-n+i) \sum^{r-1}_{j=1} p^{(e,I)}_{i,\gamma-i,j} q^{(e,I)}_{n-i,m-\gamma-n+i,r-j} +  C^{(e)}_3 \sum^{r-1}_{j=0} q^{(e,I)}_{i,\gamma-i,j} q^{(e,I)}_{n-i,m-\gamma-n+i-1,r-j-1}  \\
        & + (n-i+1) C^{(e)}_2 \sum^{r-1}_{j=1} p^{(e,I)}_{i,\gamma-i,j}p^{(e,I)}_{n-i+1,m-\gamma-n+i-1,r-j} \Bigg\}
        \end{aligned}
    \end{cases} \\
    0\leq n \leq m-1,\; for\;\;\; 1 \leq m \leq M
\end{gathered}
\end{equation}
\end{figure*}

To make the linear equations balanced (i.e. the number of unknowns equals the number of equations and the equations has unique solution), we need $N_{eqn}=N_{unk}$, which can be met by setting $0 < M - M^{(e)}_x = M^{(\nu)}_B = M^{(\nu)}_J = M^{(e)}_s < M$ for all the $e \in \mathcal{E}$ and $\nu \in \mathcal{V}$. Then, the unknowns regarding $s^0$ can be uniquely solved by $\mathbf{c} = \mathbf{A}^{-1}\mathbf{b}$.

\textit{\textbf{Remark}}:
Although there exists flexibility in the selection of $\Delta x_k$ and $\Delta t_k$, in this paper, we only consider uniform distribution of points i.e., $\Delta x_k$ and $\Delta t_k$ of all the cells are selected such that $\Delta t_{k+1} - \Delta t_{k} = \Delta t_1$ for all the $1<k<M_B$ (same for $M_J$ and $M_s$), and $\Delta x_{k+1} - \Delta x_{k} = \Delta x_1$ for all the $1<k<M_x$.

\subsection{Solution of unknowns regarding $s^r$, $r > 0$}

After solving the simplified problem at $s=0$ and identifying the unknowns regarding $s^0$, the original problem is considered which consists of the constraints \eqref{eq_ebd_PDE}, \eqref{eq_nor_ini}, \eqref{eq_nor_bd_P}, \eqref{eq_nor_bd_Q}, \eqref{eq_nor_bd_J} and \eqref{eq_nor_seam}. Similar to solving the simplified problem, by substituting solution \eqref{eq_SAS_Sol} into the constraints, a set of equations regarding the unknowns, which are formally linear, will be obtained. The main difference is that the linear equations are obtained recursively; that is, given that the unknowns regarding $s^0$, $s^1$, ..., $s^{r-1}$ are known, a set of linear equations regarding only the unknowns regarding $s^r$ are obtained, and then, the unknowns regarding $s^r$ are solved. 

Two different schemes for solving the problem are provided in the following subsections. The main difference is the modification on the PDE constraints \eqref{eq_ebd_PDE}. The first method, referred to as SAS-1, solves \eqref{eq_ebd_PDE} . While the second method, referred to as SAS-2, adopts simplification on \eqref{eq_ebd_PDE} so that the computation burden is greatly reduced.

\subsubsection{SAS-1}

Since the linear equations regarding the unknowns can be similarly obtained as in the simplified problem, the details of the treatment of the constraints are ignored. Only the linear equations regarding the unknowns regarding $s^r$ are given, assuming the unknowns regarding $s^0$, $s^1$, ..., $s^{r-1}$ are known. 

A set of linear equations \eqref{eq_SAS_PDE_coef_sr}, \eqref{eq_SAS_ini_coef_sr}, \eqref{eq_SAS_bd_P_Q_coef_sr}, \eqref{eq_SAS_bd_J_coef_sr} and \eqref{eq_SAS_seam_coef_sr} are derived. Specifically, we derive \eqref{eq_SAS_PDE_coef_sr} from the PDE \eqref{eq_ebd_PDE}, \eqref{eq_SAS_ini_coef_sr} from initial value \eqref{eq_nor_ini}, \eqref{eq_SAS_bd_P_Q_coef_sr} from boundary condition regarding $\mathcal{V_P}$ \eqref{eq_nor_bd_P} and $\mathcal{V_Q}$ \eqref{eq_nor_bd_Q}, \eqref{eq_SAS_bd_J_coef_sr} from the boundary condition regarding $\mathcal{V_J}$ \eqref{eq_nor_bd_J}, and \eqref{eq_SAS_seam_coef_sr} from seamless condition \eqref{eq_nor_seam}. 

\begin{equation} \label{eq_SAS_ini_coef_sr}
    \begin{cases}
        \sum^{M}_{m=0} \Delta x_k^m p^{(e,I)}_{m,0,r} = 0\\
        \sum^{M}_{m=0} \Delta x_k^m q^{(e,I)}_{m,0,r} = 0  
    \end{cases} \;\;\;\;\;\;\;\;\;\;\;\;\;1 \leq k \leq M^{(e)}_x
\end{equation}

\begin{equation} \label{eq_SAS_bd_P_Q_coef_sr}
\begin{gathered}
    \begin{cases}
        \sum^{M}_{m=0} \Delta t_k^{m} p^{(e,1)}_{0,m,r} = 0, \;\;\;\;\;\;\;\;\;\;\;\;\;\;\;\;\; \text{if $\nu \in \mathcal{V_P}$} \\
        \sum^{M}_{m=0}\sum^{m}_{n = 0} \Delta t_k^{m-n} q^{(e,N^{(e)})}_{n,m-n,r}  = 0,\;\;\;\;\text{if $\nu \in \mathcal{V_Q}$}
    \end{cases} \\
    1 \leq k \leq M^{(\nu)}_B
\end{gathered}
\end{equation}

\begin{equation} \label{eq_SAS_bd_J_coef_sr}
\begin{gathered}
    \begin{cases}
        \sum^{M}_{m=0} \Delta t_k^m  p^{(e^{in,1},1)}_{0,m,r} = \sum^{M}_{m=0} \Delta t_k^m  p^{(e^{in,i},1)}_{0,m,r},\;\;\;\;\;\;\;\;\;i \neq 1 \\
        \sum^{M}_{m=0} \Delta t_k^m  p^{(e^{in,1},1)}_{0,m,r} = \sum^{M}_{m=0} \sum^{m}_{n=0} \Delta t_k^{m-n} p^{\left( e^{out,j},N^{(e^{out,j})}\right)}_{n,m-n,r}  \\
        \sum_{e^{out,j}} \left[   \sum^{M}_{m=0} \sum^{m}_{n=0} \Delta t_k^{m-n} q^{\left(e^{out,j},N^{(e^{out,j})}\right)}_{n,m-n,r}   \right] - \\
        \;\;\;\;\;\;\;\;\;\;\;\;\;\;\;\;\;\;\;\;\;\;\;\;\;\;\;\;\;\;\sum_{e^{in,i}} \left[  \sum^{M}_{m=0} \Delta t_k^m  q^{(e^{in,i},1)}_{0,m,r} \right] = 0 
    \end{cases}\\
    1 \leq k \leq M^{(\nu)}_J
\end{gathered}
\end{equation}

\begin{equation} \label{eq_SAS_seam_coef_sr}
\begin{gathered} 
    \begin{cases}
        \sum^{M}_{m=0}\sum^{m}_{n = 0} \Delta t_k^{m-n} p^{(e,I)}_{n,m-n,r} = \sum^{M}_{m=0} \Delta t_k^{m} p^{\left( e,(I+1) \right)}_{0,m,r}  \\
        \sum^{M}_{m=0}\sum^{m}_{n = 0} \Delta t_k^{m-n} q^{(e,I)}_{n,m-n,r} = \sum^{M}_{m=0} \Delta t_k^{m} q^{\left( e,(I+1) \right)}_{0,m,r}
    \end{cases} \\ 
1 \leq k \leq M^{(e)}_s 
\end{gathered}    
\end{equation}

All the linear equations for identifying the unknowns regarding $s^r$ can be organized in the form of \eqref{eq_lin_sys_sr}, where $\mathbf{c}^{(r)}$ is the vector of all the unknowns, $\mathbf{A}^{(r-1)}$ is the resulting matrix and $\mathbf{b}^{(r-1)}$ is a vector of constants. 

\begin{equation} \label{eq_lin_sys_sr}
    \mathbf{A}^{(r-1)}\mathbf{c}^{(r)} = \mathbf{b}^{(r-1)}
\end{equation}

The total number of equations and unknowns can also be computed by \eqref{eq_eqn_num} and \eqref{eq_unk_num}, respectively. Select $0 < M - M^{(e)}_x = M^{(\nu)}_B = M^{(\nu)}_J = M^{(e)}_s < M$ to ensure balanced linear equations, i.e. $N_{eqn}=N_{unk}$. Then, the unknowns regarding $s^r$ can be uniquely solved by $\mathbf{c}^{(r)} = \left(\mathbf{A}^{(r-1)}\right)^{-1}\mathbf{b}^{(r-1)}$. 

\textbf{\textit{Remarks:}}

\textit{a}) It is noticeable that $\mathbf{A}^{(r)} = \mathbf{A}^{(1)}$ for any $r>1$. Hence, the inversion of  $\mathbf{A}^{(r)}$ is only computed once at $r=1$ which greatly reduces the computation burden. Moreover, when $\left\lVert \mathbf{B}^{(r-1)} \right\rVert_{infty} < \varepsilon_{\mathbf{b}}$ is detected at $r = \tilde r$, the identification of $\mathbf{c}^{(r)}$ can be terminated in advance and set $\mathbf{c}^{(r)} = \mathbf{0}$ for $ \tilde r \leq r \leq R$, which further reduces the computation burden.

\textit{b}) Eq. \eqref{eq_lin_sys} can be viewed as a special case of \eqref{eq_lin_sys_sr} at $r = 0$. Hence, SAS-1 is generalized to recursively solve \eqref{eq_lin_sys_sr} from $r = 0\rightarrow1\rightarrow2\rightarrow...\rightarrow R$. 

\subsubsection{SAS-2}

The PDE-related linear equation \eqref{eq_SAS_PDE_coef_sr} could be computationally intensive. Instead of directly using the constraint \eqref{eq_ebd_PDE} to derive \eqref{eq_SAS_PDE_coef_sr}, we linearize the friction term in segments and thus the simplified PDE is used as shown in \eqref{eq_ebd_PDE2}:

\begin{equation} \label{eq_ebd_PDE2}
    \begin{cases}
        \partial_{\Delta t} p^{(e,I)} + C^{(e)}_1 \partial_{\Delta x} q^{(e,I)} = 0 \\
        \partial_{\Delta t} q^{(e,I)} + C^{(e)}_2 \partial_{\Delta x} p^{(e,I)} + s C^{(e)}_3 C^{(e,I)}_4 q^{(e,I)} = 0
    \end{cases}
\end{equation}

\begin{equation} \label{eq_C4}
    C^{(e,I)}_4 = \frac{ q^{(e,I)}(0,0) + q^{(e,(I+1))}(0,0)}{ p^{(e,I)}(0,0) + p^{(e,(I+1))}(0,0) } 
\end{equation}

By substituting \eqref{eq_SAS_Sol} into \eqref{eq_ebd_PDE2}, the PDE-related linear equation becomes \eqref{eq_SAS_PDE_coef_sr2}, which has a much simpler expression than \eqref{eq_SAS_PDE_coef_sr} and the computation burden is also smaller.

\begin{equation} \label{eq_SAS_PDE_coef_sr2}
\begin{gathered}
    \begin{cases}
        (m-n) p^{(e,I)}_{n, m-n, r} + (n+1) C^{(e)}_1 q^{(e,I)}_{n+1, m-n-1, r} = 0 \\
        \begin{aligned}
        (m-n) q^{(e,I)}_{n, m-n, r} + (n+1) C^{(e)}_2 &p^{(e,I)}_{n+1, m-n-1, r} = \\
        &-C^{(e)}_3 C^{(e,I)}_4 q^{(e,I)}_{n, m-n-1, r-1}     
        \end{aligned}
    \end{cases} \\
    0\leq n \leq m-1,\; for\;\;\; 1 \leq m \leq M
\end{gathered}
\end{equation}

Similar to SAS-1, all the linear equations \eqref{eq_SAS_PDE_coef_sr2}, \eqref{eq_SAS_ini_coef_sr}, \eqref{eq_SAS_bd_P_Q_coef_sr}, \eqref{eq_SAS_bd_J_coef_sr} and \eqref{eq_SAS_seam_coef_sr} can be organized linear equations like \eqref{eq_lin_sys_sr} and the unknowns $\mathbf{c}^{(r)}$ can be solved.

\textbf{\textit{Remarks:}}

\textit{a}) For SAS-2, \eqref{eq_lin_sys} can also be viewed as a special case of \eqref{eq_lin_sys_sr} at $r = 0$, and SAS-2 can also be generalized to recursively solve \eqref{eq_lin_sys_sr} from $r = 0\rightarrow1\rightarrow2\rightarrow...\rightarrow R$. 

\textit{b}) $\mathbf{A}^{(r)} = \mathbf{A}^{(0)}$ for any $r>0$. Moreover, $\mathbf{A}^{(0)}$ will stay unchanged throughout the time-stepping simulation. Therefore, the inversion of $\mathbf{A}^{(0)}$ only needs to be computed once at the beginning of the simulation, which significantly reduces the computation burden.

\section{Overall simulation procedure} \label{sec:4}

\subsection{Simulation procedure for SAS}

\begin{figure}[thpb]
    \centering
    \includegraphics[scale=0.44]{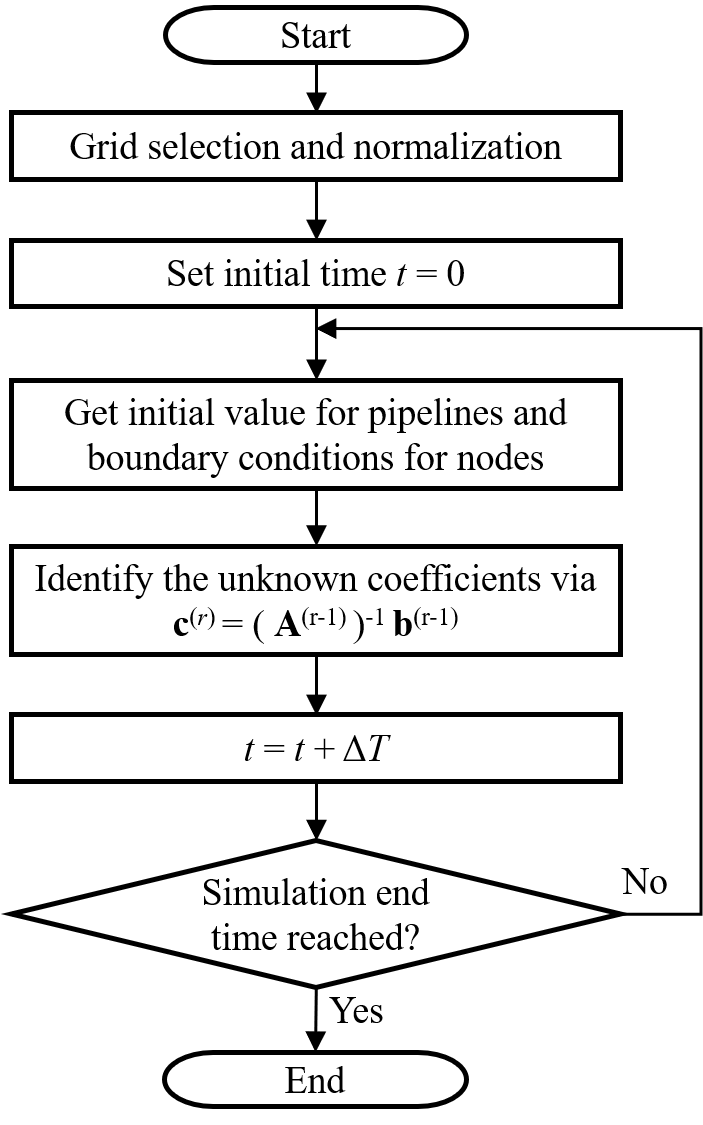}
    \caption{Simulation procedure for SAS.}
    \label{fig:SAS_sim}
\end{figure}

With the SAS approach introduced in the previous section, the key steps of the whole simulation procedure are illustrated in Fig. \ref{fig:SAS_sim}. 

In either SAS-1 or SAS-2, solving the BVP problem has been transformed to solving a set of linear equations, which is convenient by using the ubiquitous linear solvers like LU decomposition. Note that the SAS approach can consider the constraints of PDEs for the whole spatio-temporal space, while FDM only consider the same constraints at limited grid points. Thus, SAS is expected to have better accuracy and efficiency.

\subsection{Simulation procedure for implicit finite difference method}

\begin{figure}[thpb]
    \centering
    \includegraphics[scale=0.5]{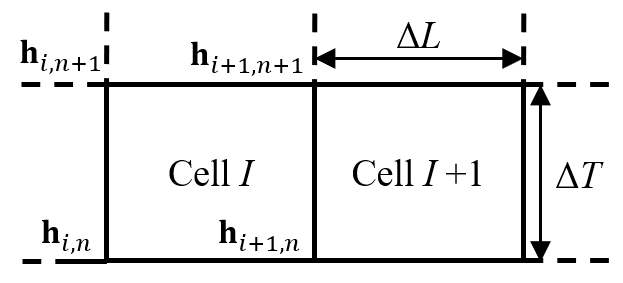}
    \caption{Grid cell of implicit FDM.}
    \label{fig:FDM_grid}
\end{figure}

The implicit FDM is used for comparison purpose and the simulation procedure is also in a ``time-stepping'' manner. In this paper we implement cell-centered method \cite{helgaker2014transient, helgaker2013modeling}, which is a widely used scheme of implicit FDM. As shown in Fig. \ref{fig:FDM_grid}, after the region of interests is discretized into a grid, the values at the four corners of each cell are used to approximate the partial derivatives, where $\mathbf{h}_{i,n}$ denotes the value of $\mathbf{h}$ at spatial index $i$ and the temporal index $n$. The partial derivatives of $\mathbf{h}$ regarding the Cell $I$ are approximated by \eqref{eq_h_deriv_appro}, which shows that it has second-order accuracy in $t$ and first-order accuracy in $x$.

\begin{equation} \label{eq_h_deriv_appro}
    \begin{cases}
        \partial_{\Delta t} \mathbf{h}^{(e,I)} \approx \frac{\mathbf{h}^{^{(e,I)}}_{i+1,n+1} + \mathbf{h}^{(e,I)}_{i,n+1} - \mathbf{h}^{(e,I)}_{i+1,n} - \mathbf{h}^{(e,I)}_{i,n} }{2 \Delta T} \\
        \partial_{\Delta x} \mathbf{h}^{(e,I)} \approx  \frac{\mathbf{h}^{(e,I)}_{i+1,n+1} - \mathbf{h}^{(e,I)}_{i,n+1} }{\Delta L}
    \end{cases} 
\end{equation}

To ensure a numerically-stable solution, the nonlinear friction term in \eqref{eq_nor_PDE} is discretized as \eqref{eq_non_term_FDM} \cite{helgaker2013modeling}.

\begin{equation}\label{eq_non_term_FDM} \small
    \frac{ q^{(e,I)} |q^{(e,I)}| }{p^{(e,I)}} \approx \frac{q^{(e,I)}_{i+1,n} + q^{(e,I)}_{i,n}}{p^{(e,I)}_{i+1,n} + p^{(e,I)}_{i,n}} \left( q^{(e,I)}_{i+1,n+1} + q^{(e,I)}_{i,n+1} - \frac{q^{(e,I)}_{i+1,n} + q^{(e,I)}_{i,n}}{2} \right)
\end{equation}

Within one time step, a set of balanced linear equations can be derived by discretizing the constraints of the normalized BVP problem via \eqref{eq_h_deriv_appro} and \eqref{eq_non_term_FDM} for all the cells in a row. Note that all the $\mathbf{h}$ at the time index $n$ are given by the initial values. Among all the $\mathbf{h}$ at the time index $n+1$, some could be given according to the boundary conditions, while the rest are unknowns to be solved from the balanced linear equations. For the gas networks, the balanced linear equations are of dimension $N_{FDM}$ as computed by:
\begin{equation} \label{eq_FDM_num}
    N_{FDM} =  \sum_{e \in \mathcal{E}} \left( 2N^{(e)} \right) 
\end{equation}

\section{Case Studies}\label{sec:5}

\subsection{Single pipeline}

The single pipeline case only has one supplying node at the inlet and one demanding node at the outlet. The configuration of the simulation is set as follows. Set $L = 2000$ m, $d = 1.016$ m, $v = 380$ m/s, $S = 0.8107$ m$^2$, $\lambda = 0.0075$, $p_b = 1\times10^6 $ Pa, and $q_b = 2\times 10^3$ kg/s. The simulation duration is $200$ s. The boundary conditions are set as \eqref{eq_case1_bd}, where the mass flow $q$ at the outlet is a sinusoidal function in time domain while the pressure $p$ at the inlet remains constant. The initial value corresponds to a stable operating condition that is compatible with the boundary conditions at $t=0$. Both SAS and implicit FDM approaches are implemented in MATLAB.

\begin{equation} \label{eq_case1_bd}
\begin{cases}
    P_B^{inlet} (t) = 6\times10^6\;\mathrm{(Pa)} \\
    Q_B^{outlet} (t) = 270 + 30 \cos{\big( \frac{\pi}{10} t \big) } \mathrm{(kg/s)}
\end{cases}
\end{equation} 

The performance of the SAS-1, SAS-2 and FDM approaches is tested in terms of accuracy and simulation speed.  For the SAS approach, $M$ is selected from $2$, $3$ and $4$ with $M_x = 1$. The result from the implicit FDM with a very tiny cell size $ [ \Delta L,\text{ }\Delta T] = [1 \text{ m},\text{ }0.001\text{ s}]$ is used as a reference for the comparison of accuracy. The metric of simulation error is defined in \eqref{eq_q_error}, where $q^{sim}_{inlet}(t)$ is simulated result $q$ at the inlet while $q^{ref}_{inlet}(t)$ is the corresponding reference result. $k$ indicates each time step. The simulation error compared with reference $ERR$ is defined as follows:

\begin{equation} \label{eq_q_error}
ERR =  \max_k \left( \frac{ \lvert q^{sim}_{inlet}(k \Delta T) - q^{ref}_{inlet}(k \Delta T)  \rvert  }{ q_B }  \right)
\end{equation}

Different cell sizes are considered for comparison. Fix $\Delta L = 400$ m, and select $\Delta T$ from $0.05$ s, $0.1$ s, $0.2$ s, $0.5 $ s and $1$ s. The simulation error of SAS-1, SAS-2, and FDM are given in Table \ref{table_case1_err_SAS-1}, Table \ref{table_case1_err_SAS-2} and Table \ref{table_case1_err_FDM}, respectively. The time consumption of SAS-1, SAS-2, and FDM are given in Table \ref{table_case1_tc_SAS-1}, Table \ref{table_case1_tc_SAS-2} and Table \ref{table_case1_tc_FDM}, respectively. The red crossing denotes that the simulation diverges. Simulation tends to diverge for the SAS approaches under $M=3$ which will be discussed and improved later.

\begin{figure}[htb]
    \centering
    \includegraphics[scale=0.8]{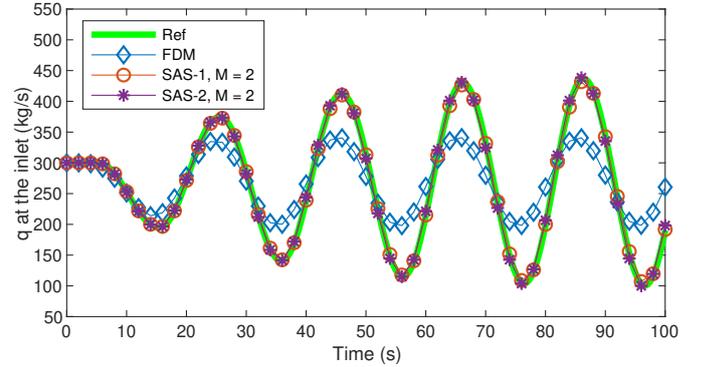}
    \caption{Simulation result of $q$ at the inlet.}
    \label{figs:case1_q_comp}
\end{figure}

The results show that SAS-1 and SAS-2 under $M=2$ and 4 are more accurate than FDM. When the cell size is smaller, the accuracy of SAS-2 is close to that of SAS-1; otherwise, it will become lower than SAS-1 but still higher than FDM. The accuracy can also be compared by visualizing the mass flow $q$ at the inlet, as shown in Fig. \ref{figs:case1_q_comp} where the cell size is $[400\text{ m}, 1\text{ s}]$ for all the approaches and the result is shown for the early 100 s. The comparison verifies the advantages of SAS approach over FDM in terms of accuracy when the same cell size is used.

\begin{table}[htb]
\small
\caption{Simulation error of SAS-1: single pipeline}
\label{table_case1_err_SAS-1}
\begin{center}
\begin{tabular}{|c||c|c|c|c|c|}
\hline
\multirow{2}{*}{\textbf{SAS-1}} & \multicolumn{5}{c|}{$log_{10}(ERR)$ \textbf{under different} $\Delta T$ } \\ \cline{2-6}  & 0.05 s &   0.1 s & 0.2 s  &  0.5 s  &  1 s  \\ 
\hhline{|=||=|=|=|=|=|}
 $M=2$ & -2.127 & -2.134 & -2.155 & -2.288 &  -2.658  \\ \cline{2-6}
\hline
 $M=3$ & \textcolor{red}{$\times$} & \textcolor{red}{$\times$}  & \textcolor{red}{$\times$}  & \textcolor{red}{$\times$}  & \textcolor{red}{$\times$}    \\ \cline{2-6} 
 \hline
 $M=4$ & -3.841 & -3.836 & -3.829 & -3.862 & -3.785   \\ \hline
\end{tabular}
\end{center}
\end{table}
\begin{table}[htb]
\small
\caption{Simulation error of SAS-2: single pipeline}
\label{table_case1_err_SAS-2}
\begin{center}
\begin{tabular}{|c||c|c|c|c|c|}
\hline
\multirow{2}{*}{\textbf{SAS-2}} & \multicolumn{5}{c|}{$log_{10}(ERR)$ \textbf{under different} $\Delta T$ } \\ \cline{2-6}  & 0.05 s &   0.1 s & 0.2 s  &  0.5 s  &  1 s  \\ 
\hhline{|=||=|=|=|=|=|}
 $M=2$ & -2.105 & -2.093 & -2.076 & -2.070 &  -2.221  \\ \cline{2-6}
\hline
 $M=3$ & -0.987 & -0.986  & -0.8241  & 0.367  & \textcolor{red}{$\times$}    \\ \cline{2-6} 
 \hline
 $M=4$ & -3.281 & -3.175 & -2.936 & -2.534 & -2.222   \\ \hline
\end{tabular}
\end{center}
\end{table}
\begin{table}[htb]
\small
\caption{Simulation error of FDM: single pipeline}
\label{table_case1_err_FDM}
\begin{center}
\begin{tabular}{|c||c|c|c|c|c|}
\hline
\multirow{2}{*}{\textbf{Method}} & \multicolumn{5}{c|}{$log_{10}(ERR)$ \textbf{under different} $\Delta T$ } \\ \cline{2-6}  & 0.05 s &   0.1 s & 0.2 s  &  0.5 s  &  1 s  \\ 
\hhline{|=||=|=|=|=|=|}
 $FDM$ & -2.089 & -1.931 & -1.705 & -1.427 &  -1.276   \\ \hline
\end{tabular}
\end{center}
\end{table}
\begin{table}[htb]
\small
\caption{Time consumption of SAS-1: single pipeline}
\label{table_case1_tc_SAS-1}
\begin{center}
\begin{tabular}{|c||c|c|c|c|c|}
\hline
\multirow{2}{*}{\textbf{SAS-1}} & \multicolumn{5}{c|}{\textbf{Time consumption under different} $\Delta T$ } \\ \cline{2-6}  & 0.05 s &   0.1 s & 0.2 s  &  0.5 s  &  1 s  \\ 
\hhline{|=||=|=|=|=|=|}
 $M=2$ & 3.634 s & 1.951 s & 0.912 s & 0.359 s &  0.180 s  \\ \cline{2-6}
\hline
 $M=3$ & \textcolor{red}{$\times$} & \textcolor{red}{$\times$}  & \textcolor{red}{$\times$}  & \textcolor{red}{$\times$}  & \textcolor{red}{$\times$}    \\ \cline{2-6} 
 \hline
 $M=4$ & 12.634 s & 6.592 s & 3.584 s & 1.332 s & 0.744 s   \\ \hline
\end{tabular}
\end{center}
\end{table}
\begin{table}[htb]
\small
\caption{Time consumption of SAS-2: single pipeline}
\label{table_case1_tc_SAS-2}
\begin{center}
\begin{tabular}{|c||c|c|c|c|c|}
\hline
\multirow{2}{*}{\textbf{SAS-2}} & \multicolumn{5}{c|}{\textbf{Time consumption under different} $\Delta T$ } \\ \cline{2-6}  & 0.05 s &   0.1 s & 0.2 s  &  0.5 s  &  1 s  \\ 
\hhline{|=||=|=|=|=|=|}
 $M=2$ & 1.367 s & 0.680 s & 0.333 s & 0.147 s &  0.069 s  \\ \cline{2-6}
\hline
 $M=3$ & 1.789 s & 0.879 s  & 0.436 s  & 0.177 s  & \textcolor{red}{$\times$}    \\ \cline{2-6} 
 \hline
 $M=4$ & 2.311 s & 1.145 s & 0.569 s & 0.236 s & 0.116 s   \\ \hline
\end{tabular}
\end{center}
\end{table}
\begin{table}[htb]
\small
\caption{Time consumption of FDM: single pipeline}
\label{table_case1_tc_FDM}
\begin{center}
\begin{tabular}{|c||c|c|c|c|c|}
\hline
\multirow{2}{*}{\textbf{Method}} & \multicolumn{5}{c|}{\textbf{Time consumption under different} $\Delta T$ } \\ \cline{2-6}  & 0.05 s &   0.1 s & 0.2 s  &  0.5 s  &  1 s  \\ 
\hhline{|=||=|=|=|=|=|}
 $FDM$ & 0.185 s & 0.074 s & 0.039 s & 0.016 s &  0.008 s   \\ \hline
\end{tabular}
\end{center}
\end{table}
\begin{table}[htb]
\small
\caption{Simulation error of SAS, $M=3$, $M_x = 2$: single pipeline}
\label{table_case1_err_M3Mx2}
\begin{center}
\begin{tabular}{|c||c|c|c|c|c|}
\hline
\multirow{2}{*}{\textbf{Method}} & \multicolumn{5}{c|}{$log_{10}(ERR)$ \textbf{under different} $\Delta T$ } \\ \cline{2-6}  & 0.05 s &   0.1 s & 0.2 s  &  0.5 s  &  1 s  \\ 
\hhline{|=||=|=|=|=|=|}
SAS-1 & -3.871 & -3.690 & -3.675 & -3.552 &  -3.219   \\ 
\hline
SAS-2 & -3.066 & -3.032 & -2.878 & -2.536 & -2.230 \\ 
 \hline
\end{tabular}
\end{center}
\end{table}

\begin{table}[htb]
\small
\caption{Time consumption of SAS, $M=3$, $M_x = 2$: single pipeline}
\label{table_case1_tc_M3Mx2}
\begin{center}
\begin{tabular}{|c||c|c|c|c|c|}
\hline
\multirow{2}{*}{\textbf{Method}} & \multicolumn{5}{c|}{\textbf{Time consumption under different} $\Delta T$ } \\ \cline{2-6}  & 0.05 s &   0.1 s & 0.2 s  &  0.5 s  &  1 s  \\ 
\hhline{|=||=|=|=|=|=|}
SAS-1 & 7.093 s & 3.882 s & 1.758 s & 0.700 s & 0.381 s   \\ 
\hline
SAS-2 &  1.809 s & 0.889 s & 0.438 s & 0.192 s & 0.090 s\\ 
 \hline
\end{tabular}
\end{center}
\end{table}

In terms of time consumption, the time consumption of FDM is generally small under the same cell size, since FDM solve a lower dimensional balanced linear equations that SAS approach within each time step. However, to reach the same degree of accuracy, FDM needs a much smaller cell than SAS which actually leads to more computation time. It can be revealed by comparing the time consumption of SAS-2 with $M=2$ and $\Delta T = 1$ s and FDM with $\Delta T = 0.05$ s. The former only needs 0.069 s to reach a higher degree of accuracy and the later needs 0.185 s. Hence, the SAS approach can admit a coarser cell  to reach higher accuracy while use shorter time consumption than FDM.

The simulation tends to diverge or reach a less accurate result for SAS with $M=3$ and $M_x = 1$, which indicates that the aforementioned efficiency may not be gained under improper configuration of SAS. One attempt is to formulate an optimization problem to minimize the objective like $\| \mathbf{b}^{(r-1)}-\mathbf{A}^{(r-1)}\mathbf{c}^{(r)}\|_2$, but it could be time-consuming. Another attempt to improve the result is to change $M_x$ from 1 to 2. The rationale behind the change of $M_x$ is that the condition number of $\mathbf{A}^{(0)}$ will decrease. For instance, when $\Delta T =$ 0.2 s, the conditional number of matrix $\mathbf{A}^{(0)}$ of SAS-2 will decrease from 979.3 to 275.2, which make the problem less ill-conditioned. In terms of the problem formulation, this change is equivalent to select only one point at the lower border for the initial value constraints in Fig. \ref{fig:ini}, and select two points at each of the left/right borders for the boundary condition constraints in Fig. \ref{fig:bon_P_Q} and the seamless condition constraints in Fig. \ref{fig:seam}. 
%Note that $M_x$ actually determines the number of points to be selected on the borders of the cell, as implied by Fig. \ref{fig:ini}, Fig. \ref{fig:bon_P_Q}, Fig. \ref{fig:seam} and $M - M^{(e)}_x = M^{(\nu)}_B = M^{(\nu)}_J = M^{(e)}_s$. When $M = 3$ and $M_x = 1$, To improve the result, it is attempted to change $M_x = 1$ to $M_x = 2$, and consequently, there are two points at the lower border in Fig. \ref{fig:ini} and one point at each of the left/right borders in Fig. \ref{fig:bon_P_Q} and Fig. \ref{fig:seam}. 
The resulting simulation error and time consumption of the SAS approaches are given in Table \ref{table_case1_err_M3Mx2} and Table \ref{table_case1_tc_M3Mx2}, respectively. The comparison shows that the accuracy is improved with $M_x = 2$. which indicates that proper selection of the points for the initial value, boundary condition and seamless condition constrains can improve the accuracy of the SAS approaches. A more comprehensive investigation on proper SAS configuration will be conducted in the future, considering not just $M_x$ but also other parameters like the cell sizes.

In general, with proper parameter selections, the SAS approaches can benefit from a coarser cell size to reach the same degree of accuracy with less time consumption. Besides, compared to SAS-1, SAS-2 can be viewed as a better approach in terms of the trade-off between accuracy and time consumption, i.e. SAS-2 uses much less computation time than SAS-1 without greatly reducing the accuracy. In addition, SAS-2 can have better convergence than SAS-1 by comparing their simulation results under $M=3$ and $M_x = 1$. Hence, SAS-2 is preferred compared to SAS-1 for the purpose of practical application.

\subsection{6-node network}

\begin{figure}[thpb]
    \centering
    \includegraphics[scale=0.4]{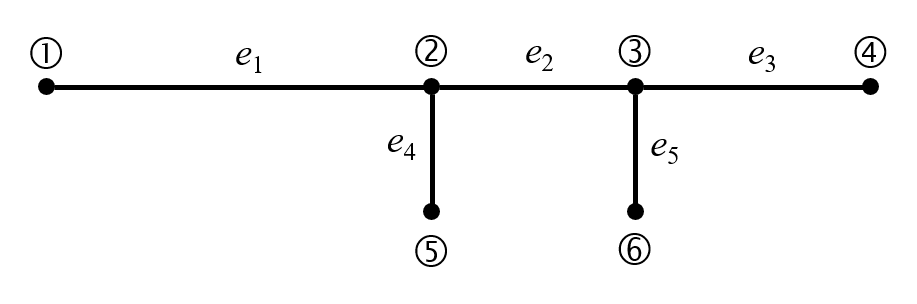}
    \caption{Diagram of 6-node network}
    \label{figure_6node}
\end{figure}

The diagram of 6-node network is shown in Fig. \ref{figure_6node}. It has five pipelines and six nodes. Supplying node include node 1, demanding nodes includes node 4, 5 and 6, and junction nodes includes node 2 and 3. Set $p_b = 6\times10^6 $ Pa, and $q_b = 2\times 10^3$ kg/s. The simulation duration is $100$ s. The parameters of the network and the grid selection is given in Table \ref{table_case2_parameter}. The boundary conditions are set as \eqref{eq_casenw_bd}. The initial value corresponds to a stable operating condition that is compatible with the boundary conditions at $t=0$.

%, where $\Delta L_{base}$ is base selection of the spatial length of the cell. Note that the change of spatial length is equivalent to scaling on the base $\Delta L_{base}$. 

\begin{table}
\caption{Parameters of Pipelines}
\label{table_case2_parameter}
\begin{center}
\begin{tabular}{|c||c|c|c|c|c|}
\hline
Pipeline & $e_1$ & $e_2$ & $e_3$ & $e_4$ & $e_5$\\
\hline
$L$ (km)  &  8 &  4 &  6 & 2 &  2 \\
\hline
$d$ (m)   & 1.2 & 1  & 1 & 1 & 1 \\
\hline 
$S$ (m$^2$)   & 1.131  &  0.785  & 0.785 & 0.785 & 0.785 \\
\hline
$\lambda$   &  0.0214 &  0.015 &   0.015 &  0.015 &  0.015 \\
\hline
$\Delta L$ (m)  & 400 & 200 & 300 & 200 & 200\\
\hline
\end{tabular}
\end{center}
\end{table}

\begin{equation} \label{eq_casenw_bd}
\begin{cases}
    P^{(1)}_B (t) = 6.96 \times10^6\;\mathrm{(Pa)} \\
    Q^{(2)}_J (t) = 0 \mathrm{(kg/s)} \\
    Q^{(3)}_J (t) = 0 \mathrm{(kg/s)} \\
    Q^{(4)}_B (t) = 125 -25 \cos{\big( \frac{\pi}{5} t \big) } \mathrm{(kg/s)} \\
    Q^{(5)}_B (t) = 210 -60 \cos{\big( \frac{\pi}{10} t \big) } \mathrm{(kg/s)} \\
    Q^{(6)}_B (t) = 125 -25 \cos{\big( \frac{\pi}{5} t \big) } \mathrm{(kg/s)}
\end{cases}
\end{equation}

The performance of SAS-1, SAS-2 and FDM is tested in terms of accuracy and simulation speed. The result from the FDM with a very tiny cell size $  \Delta L^{(e)}_{ref}=0.025 \Delta L^{(e)}, \Delta T=0.001 s$ is used as a reference for the comparison of accuracy. As inspired by the case study of single pipeline, select ($M$, $M_x$) = (2, 1), (3, 2), (4, 1) for SAS-1 and SAS-2. The accuracy is computed by applying the metric \eqref{eq_q_error} on the pipeline $e_3$.

Different temporal size $\Delta T$ are considered for comparison. Select $\Delta T$ from $0.1$ s, $0.2$ s, $0.3$ s, $0.4$ s and $0.5$ s. The simulation error of SAS-1, SAS-2, and FDM are given in Table \ref{table_case2_err_SAS-1}, Table \ref{table_case2_err_SAS-2} and Table \ref{table_case2_err_FDM}, respectively. The time consumption of SAS-1, SAS-2, and FDM are given in Table \ref{table_case2_tc_SAS-1}, Table \ref{table_case2_tc_SAS-2} and Table \ref{table_case2_tc_FDM}, respectively.

\begin{table}[thpb]
\small
\caption{Simulation error of SAS-1: 6-node network}
\label{table_case2_err_SAS-1}
\begin{center}
\begin{tabular}{|c||c|c|c|c|c|}
\hline
\multirow{2}{*}{\textbf{SAS-1}} & \multicolumn{5}{c|}{$log_{10}(ERR)$ \textbf{under different} $\Delta T$ } \\ \cline{2-6}  & 0.1 s &  0.2 s &  0.3 s & 0.4 s & 0.5 s \\ 
\hhline{|=||=|=|=|=|=|}
 $M=2$ & -2.111  &  -2.132 &  -2.166 &  -2.220 &  -2.312 \\ \cline{2-6}
\hline
 $M=3$ & -3.813  &  -3.800 &  -3.782 &  -3.725 &  -3.642 \\ \cline{2-6} 
 \hline
 $M=4$ & -3.655  &  -3.652 &  -3.705 &  -3.743 &  -3.743  \\ \hline
\end{tabular}
\end{center}
\end{table}

\begin{table}[thpb]
\small
\caption{Simulation error of SAS-2: 6-node network}
\label{table_case2_err_SAS-2}
\begin{center}
\begin{tabular}{|c||c|c|c|c|c|}
\hline
\multirow{2}{*}{\textbf{SAS-2}} & \multicolumn{5}{c|}{$log_{10}(ERR)$ \textbf{under different} $\Delta T$ } \\ \cline{2-6}  & 0.1 s &  0.2 s &  0.3 s & 0.4 s & 0.5 s \\ 
\hhline{|=||=|=|=|=|=|}
 $M=2$ & -2.096  & -2.102  &  -2.120 &  -2.150 &  -2.212 \\ \cline{2-6}
\hline
 $M=3$ & -3.255  & -3.135  &  -3.004 &  -2.906 &  -2.813 \\ \cline{2-6} 
 \hline
 $M=4$ & -3.353  & -3.173  &  -3.014 &  -2.906 &  -2.810  \\ \hline
\end{tabular}
\end{center}
\end{table}

\begin{table}[thpb]
\small
\caption{Simulation error of FDM: 6-node network}
\label{table_case2_err_FDM}
\begin{center}
\begin{tabular}{|c||c|c|c|c|c|}
\hline
\multirow{2}{*}{\textbf{Method}} & \multicolumn{5}{c|}{$log_{10}(ERR)$ \textbf{under different} $\Delta T$ } \\ \cline{2-6}  & 0.1 s &  0.2 s &  0.3 s & 0.4 s & 0.5 s \\ 
\hhline{|=||=|=|=|=|=|}
 $FDM$ & -1.943 &  -1.784 &  -1.696  &  -1.642  &  -1.602  \\ \hline
\end{tabular}
\end{center}
\end{table}

\begin{table}[thpb]
\small
\caption{Time consumption of SAS-1: 6-node network}
\label{table_case2_tc_SAS-1}
\begin{center}
\begin{tabular}{|c||c|c|c|c|c|}
\hline
\multirow{2}{*}{\textbf{SAS-1}} & \multicolumn{5}{c|}{\textbf{Time consumption under different} $\Delta T$ } \\ \cline{2-6}  & 0.1 s &  0.2 s &  0.3 s & 0.4 s & 0.5 s \\ 
\hhline{|=||=|=|=|=|=|}
 $M=2$ & 12.001 s  & 5.863 s &  3.837 s  & 2.884 s & 2.309 s  \\ \cline{2-6}
\hline
 $M=3$ & 35.562 s  & 16.878 s  & 11.251 s  & 8.405 s & 6.730 s  \\ \cline{2-6} 
 \hline
 $M=4$ & 87.118 s  & 42.911 s & 29.911 s   & 22.778 s & 18.409 s \\ \hline
\end{tabular}
\end{center}
\end{table}

\begin{table}[thpb]
\small
\caption{Time consumption of SAS-2: 6-node network}
\label{table_case2_tc_SAS-2}
\begin{center}
\begin{tabular}{|c||c|c|c|c|c|}
\hline
\multirow{2}{*}{\textbf{SAS-2}} & \multicolumn{5}{c|}{\textbf{Time consumption under different} $\Delta T$ } \\ \cline{2-6}  & 0.1 s &  0.2 s &  0.3 s & 0.4 s & 0.5 s \\ 
\hhline{|=||=|=|=|=|=|}
 $M=2$ & 4.506 s & 2.266 s  &  1.500 s & 1.173 s & 0.915 s \\ \cline{2-6}
\hline
 $M=3$ & 7.682 s & 3.842 s  &  2.543 s & 1.979 s & 1.695 s  \\ \cline{2-6} 
 \hline
 $M=4$ & 11.846 s & 6.113 s  & 4.608 s & 3.538 s & 2.785 s  \\ \hline
\end{tabular}
\end{center}
\end{table}

\begin{table}[thpb]
\small
\caption{Time consumption of FDM: 6-node network}
\label{table_case2_tc_FDM}
\begin{center}
\begin{tabular}{|c||c|c|c|c|c|}
\hline
\multirow{2}{*}{\textbf{Method}} & \multicolumn{5}{c|}{\textbf{Time consumption under different} $\Delta T$ } \\ \cline{2-6} & 0.1 s &  0.2 s &  0.3 s & 0.4 s & 0.5 s \\ 
\hhline{|=||=|=|=|=|=|}
 $FDM$ & 0.716 s &  0.379 s &   0.236 s & 0.185 s & 0.150 s  \\ \hline
\end{tabular}
\end{center}
\end{table}

\begin{figure}[thpb]
    \centering
    \includegraphics[scale=0.8]{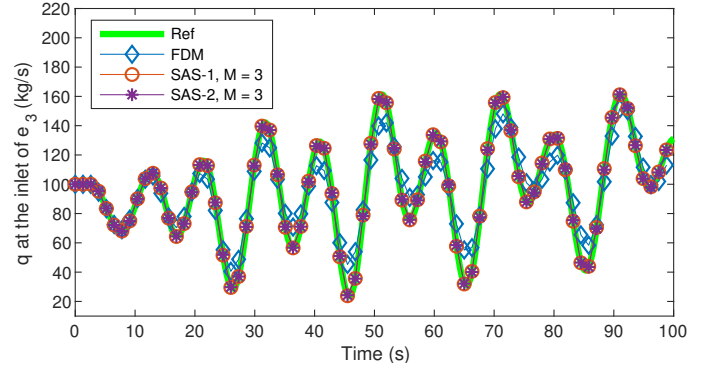}
    \caption{Simulation result of $q$ at the inlet of pipeline $e_2$.}
    \label{figs:case2_q_comp}
\end{figure}

\begin{figure}[thpb]
    \centering
    \includegraphics[scale=0.8]{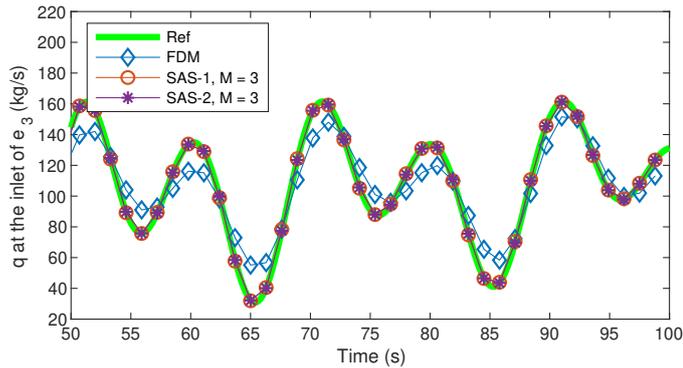}
    \caption{Simulation result of $q$ within 50 s $\leq t \leq$ 100 s.}
    \label{figs:case2_q_comp_zoomin}
\end{figure}

The comparison shows that SAS-1 and SAS-2 are more accurate than FDM, and SAS-2 is less accurate than SAS-1. The accuracy can also be compared by visualizing the mass flow $q$ at the inlet of pipeline $e_3$, as shown in Fig. \ref{figs:case2_q_comp} where the cell size is $[ \Delta L^{(e)},\text{ }0.1\text{ s}]$ for all the approaches. The result within 50 s $\leq t \leq$  100 s is also given in Fig. \ref{figs:case2_q_comp_zoomin} for better comparison. The comparison verifies the advantages of SAS approach over FDM in terms of accuracy when the same cell size is used.

In terms of time consumption, the time consumption of FDM is generally small under the same cell size at the cost of less accuracy. By reducing the cell size to $[ 0.2\Delta L^{(e)},\text{ }0.02\text{ s}]$, FDM is able to reach a high accuracy level $log_{10}(ERR) = -2.560$ with the time consumption to be 13.877 s. However, the SAS approaches can reach even higher accuracy level with less time consumption. When the cell size is $[ \Delta L^{(e)},\text{ }0.5\text{ s}]$,  SAS-1 with $M=3$ and $M_x = 2$ can reach $log_{10}(ERR) = -3.642$ within 6.730 s. When the cell size is $[ \Delta L^{(e)},\text{ }0.2\text{ s}]$,  SAS-2 with $M=3$ and $M_x = 2$ can reach $log_{10}(ERR) = -3.135$ within 3.842 s. Hence, the SAS approach can admit a coarser cell to reach higher accuracy while use shorter time consumption than FDM.

In general, for the case of 6-node network, the SAS approaches can also reach a higher accuracy level with less time consumption compared to FDM. SAS-2 can greatly reduce the time consumption with minor impact on accuracy.

\section{Conclusion and Future Works}\label{sec:6}

A novel semi-analytical solutions (SAS) approach is proposed in this paper for the simulation of transient flow in the gas networks with two different SAS schemes are provided. Both schemes can provide high-order approximate solutions to the transient flow and thus has advantages in accuracy and efficiency. In the SAS algorithm, a set of balanced linear equations from the constraints of PDE, initial value condition, boundary condition and seamless condition is formulated and solved. The performance of the proposed approach is compared with the implicit FDM, which validates its efficiency. The SAS-2 scheme shows promising potentials industrial applications.

The future works includes but not limited to the following directions:  (\textit{i}) how to properly select the parameters like the cell size, the order $M$ and $M_x$ to ensure both efficiency and accuracy, (\textit{ii}) how to extend the proposed approach to more complex pipeline models, (\textit{iii}) adaptive grid cell division for better error control.

%% The Appendices part is started with the command \appendix;
%% appendix sections are then done as normal sections
%\appendix

%\section{Parameters of 6-node networks and grid selection}
%\label{Appendix A}

%% If you have bibdatabase file and want bibtex to generate the
%% bibitems, please use
%%
 \bibliographystyle{elsarticle-num} 
 \bibliography{cas-refs}

%% else use the following coding to input the bibitems directly in the
%% TeX file.

% \begin{thebibliography}{00}

% %% \bibitem{label}
% %% Text of bibliographic item

% \bibitem{}

% \end{thebibliography}
\end{document}